\documentclass[12pt,reqno]{amsart}

\usepackage{multirow}
\usepackage{hhline}

\usepackage{mathtools}
\usepackage{times}
\usepackage[T1]{fontenc}
\usepackage{mathrsfs}
\usepackage{latexsym}
\usepackage[dvips]{graphics}
\usepackage[titletoc, title]{appendix}
\usepackage{amsmath,amsfonts,amsthm,amssymb,amscd}
\usepackage[dvipsnames]{xcolor}
\usepackage{hyperref}
\usepackage{amsmath}
\usepackage[utf8]{inputenc}

\usepackage{color}
\usepackage{breakurl}

\usepackage{comment}
\newcommand{\bburl}[1]{\textcolor{blue}{\url{#1}}}

\newtheorem{thm}{Theorem}[section]

\newtheorem{cor}[thm]{Corollary}

\newtheorem{lem}[thm]{Lemma}

\newtheorem{que}[thm]{Question}
\newtheorem{defi}[thm]{Definition}
\newtheorem{rek}[thm]{Remark}

\newtheorem*{cond*}{Condition A}

\email{\textcolor{blue}{\href{mailto:mshiliaev@tamu.edu}{mshiliaev@tamu.edu}}}
\address{Department of Mathematics\\ Texas A\&M University, College Station, TX 77843, USA.}

\usepackage[utf8]{inputenc}

\DeclareFixedFont{\ttb}{T1}{txtt}{bx}{n}{12} 
\DeclareFixedFont{\ttm}{T1}{txtt}{m}{n}{12}  

\usepackage{color}
\definecolor{deepblue}{rgb}{0,0,0.5}
\definecolor{deepred}{rgb}{0.6,0,0}
\definecolor{deepgreen}{rgb}{0,0.5,0}

\usepackage{listings}

\newcommand\pythonstyle{\lstset{
language=Python,
basicstyle=\ttm,
morekeywords={self},              
keywordstyle=\ttb\color{deepblue},
emph={MyClass,__init__},          
emphstyle=\ttb\color{deepred},    
stringstyle=\color{deepgreen},
frame=tb,                         
showstringspaces=false
}}

\lstnewenvironment{python}[1][]
{
\pythonstyle
\lstset{#1}
}
{}

\newcommand\pythoninline[1]{{\pythonstyle\lstinline!#1!}}

\numberwithin{equation}{section}

\DeclareFontFamily{U}{mathx}{}
\DeclareFontShape{U}{mathx}{m}{n}{<-> mathx10}{}
\DeclareSymbolFont{mathx}{U}{mathx}{m}{n}
\DeclareMathAccent{\widehat}{0}{mathx}{"70}
\DeclareMathAccent{\widecheck}{0}{mathx}{"71}

\begin{document}
\subjclass[2020]{41A65, 46B15}
\keywords{Thresholding Greedy Algorithm, Schreier unconditional, Schreier Families.}
\title{On unconditionality and higher-order Schreier unconditionality}

\author{Shiliaev Mark}

\begin{abstract}
     Let $X$ be a Banach space, $(e_n)_{n=1}^\infty$ be its basis, and  $S_\alpha$ be a Schreier family of order alpha. We introduce Condition A which is a weaker version of the Continuum Hypothesis. Granted Condition A, we show that if the basis $(e_n)$ is $S_\alpha$-unconditional for every countable ordinal alpha, then it is unconditional.

\end{abstract}

\thanks{This work was partially supported by the College of Arts \& Sciences at Texas A\&M University.}

\maketitle

\keywords{}

\maketitle

\section{Introduction}
Schreier families were introduced by Alspach and Argyros \cite{argyros} as a way to measure the complexity of weakly null sequences by assigning each sequence a countable ordinal. 
Before defining the Schreier sets, we need the following notation. For any $A, B\subset\mathbb{N},$ $A<B$ means $a<b$ for all $a\in A$ and $b\in B.$
For any $A\subset \mathbb{N}$ and $b\in\mathbb{N}$, $A>b$ means $A>\{b\}.$ 

We define the Schreier family $S_\beta$ for each countable ordinal $\beta$ recursively. 
The set $S_0$ consists of all singletons and the empty set.
Assume that for some $\beta<\omega_1$, sets $S_{\alpha}$ have been defined for all $\alpha<\beta.$
If $\beta=\gamma+1$, then
$$S_{\beta} \coloneqq \Bigl\{\bigcup_{i\ =\ 1}^mE_i:m\leq E_1<E_2<\dots<E_m,E_i\in S_\gamma, \forall1\leq i\leq m\Bigr\}.$$
If $\beta$ is a limit ordinal, then fix a sequence of successor ordinals $(\beta_n)_{n=1}^\infty$, that increases to $\beta$ with $S_{\beta_n}\subset S_{\beta_{n+1}}$ for all $n\ge1$.
Define $$S_\beta\ =\ \{E:\text{for some } m\ge1, m\le E\in S_{\beta_m}\}.$$ 
The sequence $(\beta_m)_{m=1}^\infty$ is called the $\beta$-approximating sequence.

Let $X$ be an infinite-dimensional Banach space over the field $\mathbb{F}$ and $X^\ast$ be its dual space. 
A semi-normalized basis is a collection of vectors $(e_n)^\infty_{n=1}\subset X$ such that\begin{itemize}
    \item[i)]$\overline{\mbox{span}\{(e^\ast_n)^\infty_{n=1}\}}\ =\ X$,
    \item[ii)] there exists a unique sequence $(e^\ast_n)^\infty_{n=1}\subset X^\ast$ satisfying $e^\ast_n(e_m)\ =\ \delta_{n,m}$,
    \item[iii)] there exist $c_1,c_2>0$ such that $$0<c_1:\ =\ \inf_n\{ \| e_n\|, \| e_n^\ast \| \} \le \sup_n\{ \| e_n\|, \| e_n^\ast \| \} :\ =\  c_2 <\infty.$$
\end{itemize}
 Konyagin and Temlyakov \cite{Kon} studied an order-free method of approximation and introduced the notion of greedy bases. 
For any $x\in X$, a finite set $\Lambda(x)$ is said to be a greedy set of $x$ if $\min_{n\in\Lambda}|e_n^\ast(x)|\ge\max_{n\notin\Lambda(x)}|e_n^\ast(x)|$.
We use $G(x,m)$ to denote the set of all greedy sets of $x$ of order $m$.
For any $A\subset\mathbb{N},$ let $P_A(\sum_{n}e_n^\ast(x)e_n)\ =\ \sum_{n\in A}e_n^\ast(x)e_n.$
A basis $(e_n)$ of a Banach space $X$ is called greedy if there is a $C\ge1$ such that for all $x\in X$, $m\in\mathbb{N},$ and $\Lambda\in G(x,m),$
 $$||x-P_\Lambda(x)||\leq C\inf\Bigl\{ \Big\|  x - \sum_{n\in A} a_ne_n  \Big\| 
:|A| \le m, (a_n)\subset\mathbb{K} \Bigr\}.$$
On the other hand, a basis $(e_n)$ is unconditional if there exists a $C\ge1$ such that $$||x-P_A(x)||\leq C||x||, \forall x\in X, \forall A\subset\mathbb{N}.$$
Both greediness and unconditionality are strong properties, so researchers have investigated many weaker but desirable properties \cite{ArgGas, ArgGod, ArgMer, Dil, DilKh, Odel}. This paper focuses on the notions of $\mathcal{F}$-greedy and $\mathcal{F}$-unconditional. A family $\mathcal{F}$ of finite subsets is hereditary if for any $A\in \mathcal{F},$ $B\subset A$ implies $B\in \mathcal{F}$.
\begin{defi}
    Let $\mathcal{F}$ be a hereditary family of finite subsets of $\mathbb{N}$. A basis $(e_n)$ is said to be $\mathcal{F}$-greedy if there exists $C\ge1$ such that for all $x\in X$, $m\in\mathbb{N}$ and $\Lambda\in G(x,m)$
 $$||x-P_\Lambda(x)||\leq C\inf\Bigl\{ \Big\|  x - \sum_{n\in A} a_ne_n  \Big\| 
:|A| < m,A\in\mathcal{F}, (a_n)\subset\mathbb{K} \Bigr\}.$$
\end{defi}
\begin{defi}
A basis $(e_n)$ is said to be $\mathcal{F}$-unconditional if there exists $C\ge1$ such that
 $$||x-P_A(x)||\leq C||x||, \forall x\in X, \forall A\in\mathcal{F}.$$
\end{defi}

Letting $\mathcal{F}=S_\alpha$ (the Scheier family of order $\alpha$), we get $S_\alpha$-greediness and $S_\alpha$-unconditionality, whose relationship was studied in \cite{Chu_1}. 

\begin{lem}\cite[Theorem 1.11]{Chu_1}
    Fix $\alpha\in\omega_1$.
    \begin{itemize}
        \item[1)] A basis is greedy/unconditonal if and only if there is a $C\ge1$ such that the basis is $C-S_{\alpha+n}$-greedy/unconditonal for all $n\in\mathbb{N}$.
        \item[2)] There exists a basis that is $S_{\alpha+n}$-greedy/unconditonal for all $n\in\mathbb{N}$, but is not greedy/unconditonal.
    \end{itemize}
\end{lem}
A natural next step is trying to answer the following question asked in \cite{Chu_1}.
\begin{que}\label{question}
 If a basis is $S_\alpha$-greedy/unconditional for every $\alpha \in \omega_1$, then is it 
 \\
 greedy/unconditional?   
\end{que}
This paper answers Question 1.4. positively under a condition, which we will now describe.

Let $\mathcal{G}$ denote  the set of all non-decreasing functions $\mathbb{N}\xrightarrow{}\mathbb{N}$.
\begin{defi}
    A set $H\subset \mathcal{G}$ is uniformly bounded by a function $g \in \mathcal{G}$, if for every $h\in H$, there exists $N_h\in \mathbb{N},$ such that for all $n>N_h,$ we have $h(n)<g(n).$
    
    A set $H\subset \mathcal{G}$ is uniformly bounded if it is uniformly bounded by some $g\in \mathcal{G}$.
\end{defi}
We are ready to state the condition.
\begin{cond*}
    There exists a set $\mathcal{H}\subset \mathcal{G}$ of cardinality $\aleph_1$ that is not uniformly bounded.
\end{cond*}
    
    We note that Condition A is a weaker version of Continuum Hypothesis, as shown in \cite{Blass}. In our proof, we define a function that gives us insight into the structure of the ordinal tree and Schreier sets. 

If Condition A is false, then we create a non-trivial family $S_{\omega_1}$  such that the following lemma holds.
\begin{lem}\label{anti_premain}
    If a basis is $S_{\omega_1}$-greedy/unconditional, then it is $S_\beta$-greedy/unconditional for all $\beta <\omega_1$.
\end{lem}

This family $S_{\omega_1}$ is limited in the following sense.
\begin{lem}
    For any $n\in \mathbb{N},$ there exists $m\ge n$ such that $\{n,n+1,\dots m\}\notin S_{\omega_1}.$
\end{lem}
This limitation gives us the intuition that there might exist a basis that is
$S_{\omega_1}$-
greedy/unconditional, but not greedy/unconditional. By Lemma  such basis would be a counter-example to Question 1.4.

The paper is structured as follows: Section 2 defines a function that describes the structure of the Ordinal tree; Section 3 studies the relation between that function and Schreier sets; Section 4 outlines Condition A and defines the approximating sequences; Section 5 positively answers Question 1.4. granted Condition A; Section 6 provides a building stone to give a negative answer to Question 1.4. if Condition A fails.
\section{Function describing the ordinal tree}
In this section, we first define a function $F: \mathbb{N}\times\omega_1 \xrightarrow{} \mathbb{N}\cup\{0\}$, based on the Ordinal tree. Then we show some relevant properties of that function.

For all natural numbers $n$ we define an oriented graph $G_n$, whose vertices are countable ordinals.
There are two types of edges called successor edges $S$ and limit edges $L$. For every successor vertex $\alpha+1$ we draw an edge from $\alpha+1$ to $\alpha$.
For every limit vertex $\alpha$ we draw arrows from $\alpha$ to the first $n$ elements of the $\alpha$-approximating sequence.

 We define the function $F: \mathbb{N}\times\omega_1 \xrightarrow{} \mathbb{N}\cup\{0\}$ 
 as follows: $$F(n,\alpha)\ =\ \max_{\substack{\mbox{Paths P in }\\ G_n\mbox{ from }\alpha \mbox{ to } 0}}(S \cap E(P)). $$
 We will now formulate some properties of $F$ that will be used in due course.
 \begin{lem}\label{F3Prop}

     For all natural numbers $n$ and countable ordinals $\alpha$:
     \begin{enumerate}
         \item $F(n,0)=0$,
         \item $F(n,\alpha+1)=F(n,\alpha)+1$,
         \item $F(n,\alpha)=\sup_{m\leq n}F(n,\alpha_m ) \mbox{ if } \alpha \mbox{ is a limit ordinal}.$
     \end{enumerate}
     \begin{proof}
          The function $F$ has Property (1) because the length of any path from $0$ to itself is $0$.
          $F$ has Property (2), because for every path $P$ from $\alpha+1$ to $0$, $\{\alpha+1,\alpha\}\in E(P)\cap S$.
          $F$ has Property (3), because for every path $P$ from $\alpha$ to $0$, there exists $m \in \{1,\dots,n\}$, such that $\{\alpha,\alpha+1\}\in E(P)\cap L$.
     \end{proof}
 \end{lem}

\begin{cor}\label{F_more_than_m}
    Let $n,m \in \mathbb{N}$ and $\alpha \in \omega_1$.
    Then $F(n,\alpha+m)\ge m$. 
\end{cor}

\begin{lem}\label{F_ordered}
    Let $\alpha < \beta$ be two countable ordinals.
    Then there exists $N\in \mathbb{N}$ such that $F(n,\alpha)<F(n,\beta),\mbox{ for all } n>N.$
    \begin{proof}
        We prove this lemma by transfinite induction on $\beta$.
        The base case $\beta = 1$ is trivial. 
        Inductive hypothesis: suppose the lemma holds for all $\gamma<\beta$.
        
        If $\beta$ is a successor ordinal we have $\beta = \gamma +1$, hence $\alpha \leq \gamma$.
        If $a = \gamma$, then by Property (2) from Lemma \ref{F3Prop}, $$F(n,\beta)\ =\ F(n,\alpha)+1\ >\ F(n,\alpha).$$
        If $\gamma > \alpha$, then by the inductive hypothesis, there exists $N$, such that $$F(n,\alpha)\ <\ F(n,\gamma),\ \forall n>N.$$
        Since $F(n,\beta)=F(n,\gamma)+1,$ we get $$F(n,\alpha)\ <\ F(n,\gamma)\ <\ F(n,\beta),\ \forall n>N.$$
        
        If $\beta$ is a limit ordinal, we consider the $\beta$-approximating sequence $(\beta_n)_{n=1}^\infty$.
        Choose $m$ such that $\beta_m > \alpha$. Then by the inductive hypothesis there exists $N$ such that
        $$ F(n,\alpha)\ <\ F(n,\beta_m),\mbox{ }\forall n>N.$$
        By Property (3) of $F$, $$F(n,\beta_m)\ \leq\ F(n,\beta),\mbox{ }\forall n \ge m.$$
        Therefore, $$  F(n,\alpha)<F(n,\beta_m)\ \leq\ F(n,\beta),\mbox{ }\forall n > \max(m,N).$$
    \end{proof}
\end{lem}
\begin{cor}
    If a countable ordinal $\alpha \ge \omega$, then 
    $$\lim{_{n \xrightarrow{} \infty}}F(n,\alpha)\ \ge\ \lim{_{n \xrightarrow{} \infty}}F(n,\omega)\ =\ \infty.$$
\end{cor}

\section{Relationship between $F(n,\alpha)$ and $S_\alpha$}
In this section, we will deduce properties $S_\alpha$ from properties of $F(n,\alpha)$. In particular, we get Corollary \ref{3Main}, which we will later use to positively answer Question \ref{question}.

Recall that $\mathcal{P}_{<\infty}(A)$ denotes the set of all finite subsets of $A$.
\begin{lem}\label{unb_F}
    Let $A \subset \omega_1$ and fix $n \in \mathbb{N}$. If $\ \sup_{\alpha \in A}F(n,\alpha)=\infty,$ then $$\mathcal{P}_{<\infty}(\mathbb{N}_{\ge n}) \subset \bigcup_{\alpha \in A}S_\alpha.$$
    \begin{proof}
        All Schreier families are hereditary, so it suffices to show that
        $$\{n, n+1,... ,m\} \in \bigcup_{\alpha \in A}S_\alpha,\mbox{ } \forall m \ge n.$$ This would follow from $$ \{n, n+1,... ,n+F(n,\alpha)\}\in S_\alpha, \mbox{ } \forall \alpha\in \omega_1,$$ which trivially follows from the definition of $F$, and the simple observation: $$\{n, n+1,... ,k\}\in S_\beta, \implies \{n, n+1,... ,k+1\}\in S_{\beta+1}, \mbox{ }\forall \beta \in \omega_1.$$
    \end{proof}
\end{lem}
\begin{defi}
    A set $H\subset \mathcal{G}$ is uniformly bounded by a function $g \in \mathcal{G}$, if for every $h\in H$, there exists $N_h\in \mathbb{N},$ such that for all $n>N_h,$ we have $h(n)<g(n).$
    
    A set $H\subset \mathcal{G}$ is said to be uniformly bounded if it is uniformly bounded by some $g\in \mathcal{G}$.
    
\end{defi}
\begin{defi}
    For a subset $A\subset \omega_1$ we define a set of functions $H_A\subset\mathcal{G}$ as  $$H_A\ =\ \{F(n,\alpha)|\alpha\in A\}.$$ 
\end{defi}
\begin{lem}\label{bound_F}
    Let $A \subset \omega_1,$ with $|A|=\aleph_1$ and for every $n \in \mathbb{N},$ $\sup_{\alpha \in A}F(n,\alpha)<\infty.$ Then $H_{\omega_1}$ is uniformly bounded.
    \begin{proof}
        Define $g: \mathbb{N} \xrightarrow{}\mathbb{N},$ as  $$g(n)\ =\ \sup_{\alpha \in A}F(n,\alpha).$$ Now fix $\alpha \in \omega_1$.
        Since $|A|=\aleph_1$, there exists $\beta \in A$, such that $\beta > \alpha$.
        Then by Lemma \ref{F_ordered} there exists $N \in \mathbb{N}$ such that $$F(n, \alpha)\ <\ F(n, \beta),\mbox{ }\forall n > N.$$ And by definition of $g$: $$F(n,\beta)\ \leq\ g(n),\mbox{ }\forall n\in \mathbb{N}.$$
        Hence, $$ F(n, \alpha)\ <\ F(n, \beta)\ \leq\ g(n),\mbox{ }\forall n > N.$$
        So $H_{\omega_1}$ is uniformly bounded by $g$.
    \end{proof}
\end{lem}
By combining Lemma \ref{unb_F} and Lemma \ref{bound_F} we get the following corollary.
\begin{cor}\label{3Main}
    If $H_{\omega_1}$ is not uniformly bounded, then for all $ A\subset \omega_1$, with $|A|=\aleph_1$, there exists $n \in \mathbb{N}$ such that $\mathcal{P}_{<\infty}(\mathbb{N}_{\ge n}) \subset \bigcup_{\alpha \in A}S_\alpha.$
\end{cor}
\section{Condition A and approximating sequences}
In this section, we introduce Condition A, then we define the approximating sequences for each limit ordinal granted Condition A.
\begin{lem}
    Let $\alpha,\beta \in \omega_1$ and $m,n \in \mathbb{N}$, with $m \ge n$ and $S_\alpha \subset S_\beta$, then $S_{\alpha+n}\subset S_{\beta+m}$.
    \begin{proof}
        It suffices to show that  $S_{\alpha+1} \subset S_{\beta+1}.$ Let's fix some $A\in S_{\alpha+1}.$ By definition of Schreier sets $$A\ =\ \bigcup_{i=1}^mE_i:\forall i\in\mathbb{N},E_i\in S_\alpha.$$ Since $S_\alpha \subset S_\beta$, any $E_i$ in $ S_\alpha$, is also in $S_\beta.$ Thus $A\in S_{\beta+1}$
    \end{proof}
\end{lem}
\begin{cor}\label{cor_approx}
    Fix a limit ordinal $\beta$ and $g\in\mathcal{G}.$
    If $(\beta_n)_{n=1}^\infty$ can be a $\beta$-approximating sequence, then so can $(\beta_n+g(n))_{n=1}^\infty$.
\end{cor}
The following condition trivially follows from the Continuum Hypothesis, but doesn't imply it, as shown in \cite{Blass}.
\begin{cond*}
    There exists a set $\mathcal{H}\subset \mathcal{G}$ of cardinality $\aleph_1$ that is not uniformly bounded.
\end{cond*}

The rest of Section 4 and Section 5 assume Condition A holds.
So, let's fix any such set $\mathcal{H}$. Since $|\mathcal{H}|=\aleph_1$, there exists a bijection $\varphi$ from the set of countable limit ordinals to $\mathcal{H}$.
We will denote $\varphi(\alpha)\in$ $\mathcal{H}$ as $h_\alpha$.

\begin{lem}\label{Choose_approx}
    There exists a way to choose an approximating sequence for every limit ordinal, so that for any limit ordinal $\alpha$ and $n\in\mathbb{N}$, $F(n,\alpha)\ge h_\alpha(n).$
    \begin{proof}
        We will prove this lemma by transfinite induction on limit ordinals.
        The base case is $\alpha=\omega.$ We pick the $\alpha$-approximating sequence 
        $(\alpha_n)_{n=1}^\infty$  defined term wise by $\alpha_n=h_w(n)$.
        Then by Property (3) from Lemma \ref{F3Prop}
        $$F(n,\alpha)\ \ge\ F(n,h_\alpha(n))\ =\ h_\alpha(n),\ \forall n\in \mathbb{N}.$$ Now, assume we have picked an approximating sequence for every limit ordinal
        smaller that $\alpha.$ Then there must exist some $\alpha$-approximating
        sequence, we denote this sequence by $(\alpha_n)_{n=1}^\infty.$ Define the sequence $(\alpha_n^\prime)_{n=1}^\infty$ term wise as $\alpha_n^\prime=\alpha_n+h_\alpha(n).$ By Corollary \ref{cor_approx} $(\alpha_n^\prime)_{n=1}^\infty$ can be an $\alpha$-approximating sequence. By  Corollary \ref{F_more_than_m}
        $$\ F(n,\alpha^\prime_n)\ =\ F(n,\alpha_n+h_\alpha(n))\ \ge\ h_\alpha(n),\ \forall n\in \mathbb{N}.$$
        So we pick $(\alpha_n^\prime)_{n=1}^\infty$ as the $\alpha$-approximating sequence.
    \end{proof}
\end{lem}
In the next section we answer Question 1.4 based on this construction of approximating sequences.
\section{Answer to the question 1.4}
In this section, we assume Condition A and provide a positive answer to Question \ref{question}.
We will omit the proof of the following straightforward lemma.
\begin{lem}\label{premain}
    If a basis is $\mathcal{P}_{<\infty}(\mathbb{N}_{\ge n})$-greedy/unconditional, then it is greedy/unconditional.
\end{lem}
\begin{thm}
    If a basis is $S_\alpha$-greedy/unconditional for every $\alpha \in \omega_1$, then it is greedy/unconditional. 
    \begin{proof}
        For every ordinal $\alpha$ we will use $n_\alpha$ to denote the smallest natural number for which the basis is $n_\alpha$-$S_\alpha$-greedy/unconditional.
        By the Pigeonhole principle, there exists a natural number $n$, such that the basis is $n$-$S_\alpha$-greedy/unconditional for uncountably many ordinals $\alpha$.
        Let's denote the set of those ordinals by $A$. Then the basis is $$n- \bigcup_{\alpha \in A}S_\alpha-\mbox{greedy/unconditional.}$$
        By Lemma \ref{Choose_approx} $H_{\omega_1}$ is not uniformly bounded. So by Corollary \ref{3Main}, there exists $m \in \mathbb{N}$, such that  $$\mathcal{P}_{<\infty}(\mathbb{N}_{\ge m}) \subset \bigcup_{\alpha \in A}S_\alpha.$$ 
        Hence, the basis is $n$-$ \mathcal{P}_{<\infty}(\mathbb{N}_{\ge m})$-greedy/unconditional.
        So by Lemma \ref{premain} the basis is greedy/unconditional.
    \end{proof}
\end{thm}

\section{Negation of Condition A}
This section is about the case where Condition A fails. We first define a function $G: \mathbb{N}\times\omega_1 \xrightarrow{} \mathbb{N}$, in a similar manner to the function $F$. Then we use that function $G$ to define a set $S_{\omega_1}$ and show some properties of $S_{\omega_1}$ that justify our definition.

We define the function $G: \mathbb{N}\times\omega_1 \xrightarrow{} \mathbb{N}$ 
 pointwise to satisfy the following: 
 \begin{enumerate}
     \item $[n, n+1,\dots G(n,\alpha)]\in S_\alpha,$
     \item $[n, n+1,\dots G(n,\alpha)+1]\notin S_\alpha.$
 \end{enumerate}
\begin{rek}
    
Let $n\in \mathbb{N}$, $\alpha\in \omega_1$. Then $G(n,\alpha)$ is the largest value $m$ for which $\{n,n+1,\dots,m\}\in S_\alpha.$
\end{rek}
We denote the set of functions $H_{gen}\coloneq\{ G(n,\alpha)\ |\ \alpha\in\omega_1\}.$
\begin{lem}
    $H_{gen}$ is uniformly bounded.
    \begin{proof}
        The cardinality of $H_{gen}$ is $\aleph_1$. So by negation of Condition A, $H_{gen}$ is uniformly bounded.
    \end{proof}
\end{lem}

\begin{lem}
    There exist $A\subset\omega_1$ and $g\in \mathcal{G}$,
    such that $|A|=\aleph_1$ and $G(n, \alpha)\leq g(n)$ for all $n\in \mathbb{N}$ and $\alpha\in A$.
    \begin{proof}
        Fix any $g$ that uniformly
        bounds $H_{gen}$. Then for each $\alpha\in\omega_1,$ there exists $N_\alpha\in\mathbb{N}$
        such that $$G(n,\alpha)\leq g(n),\ 
        \forall n\ge N_\alpha.$$ 
        By the Pigeonhole principle there exists  $m\in\mathbb{N}$  such that for $m=N_\alpha$ uncountably many $\alpha\in\omega_1$. Let's denote the set of those ordinals as $A$. Then $$G(n,\alpha)\ \leq\ g(n),\ \forall n\ge m,\ \forall \alpha\in A.$$ Since $G$ is an increasing function, $$G(n,\alpha)\ \leq\  G(m,\alpha)\ \leq\ g(m),\ \forall n< m,\ \forall \alpha\in A.$$ Now we define $g^\prime$ pointwise as, $$g^\prime(n)\ =\ \max(g(m),g(n)).$$ Then $$G(n,\alpha)\ \leq\ g^\prime(n),\ \forall n\in \mathbb{N},\ \forall \alpha\in A.$$
    \end{proof}
\end{lem}
\begin{defi}\label{S_omega1}
    Fix any uncountable $A\subset\omega_1$, such that there exists $g\in \mathcal{G}$ with $G(n, \alpha)\leq g(n),\ $  for all $n\in \mathbb{N}$ and $\alpha\in A$. Then define $$S_{\omega_1}\ =\ \bigcup_{\alpha\in A}S_\alpha.$$
    We call $A$ the defining set of $S_{\omega_1}$ and $g$ is called a bounding function of $S_{\omega_1}$.
\end{defi}
Note that $S_{\omega_1}$ is a union of hereditary sets, so $S_{\omega_1}$ is hereditary. Any set $S_{w_1}$ also satisfies the following lemmas.
\begin{lem}
    If a basis is $S_{\omega_1}$-greedy/unconditional, then it is $S_\alpha$-greedy/unconditional for all $\alpha <\omega_1$.
    \begin{proof}
        Fix $A$ the defining set of $S_{\omega_1}$. Since $|A|=\aleph_1,$ for all $\alpha\in\omega_1$ there exist $\beta\in A$, 
        such that $\beta>\alpha.$ So $S_{\beta}$-greedy/unconditional implies $S_\alpha$-greedy/unconditional.
        Also $S_\beta\subset S_{\omega_1}$, so $S_\beta-$greedy/unconditional implies $S_\alpha$-greedy/unconditional. Thus $S_{\omega_1}\mbox{-greedy/unconditional implies}$ $S_{\alpha}$-greedy/unconditional
        \end{proof}
\end{lem}
\begin{lem}
    For any $n\in \mathbb{N},$ there exists $m\ge n$ such that $\{n,n+1,\dots m\}\notin S_{\omega_1}.$
    \begin{proof}
        Fix $g$ a bounding function of $S_{\omega_1}$ and $A$ the defining set of $S_{\omega_1}$. For any $n\in\mathbb{N}$
        $$[n, n+1,\dots G(n,\alpha)+1]\notin S_\alpha,\  \forall\alpha\in\omega_1,$$
        and
        $$G(n,\alpha)\ \leq\  g(n),\ \forall \alpha\in A.$$
        So for any $n\in\mathbb{N}$ $$[n, n+1,\dots g(n)+1]\notin S_\alpha,\ \forall \alpha\in A.$$
        So by Definition \ref{S_omega1} for any $n\in\mathbb{N}$ $$[n, n+1,\dots g(n)+1]\notin S_{\omega_1}.$$
    \end{proof}
\end{lem}

\ \\
\end{document}